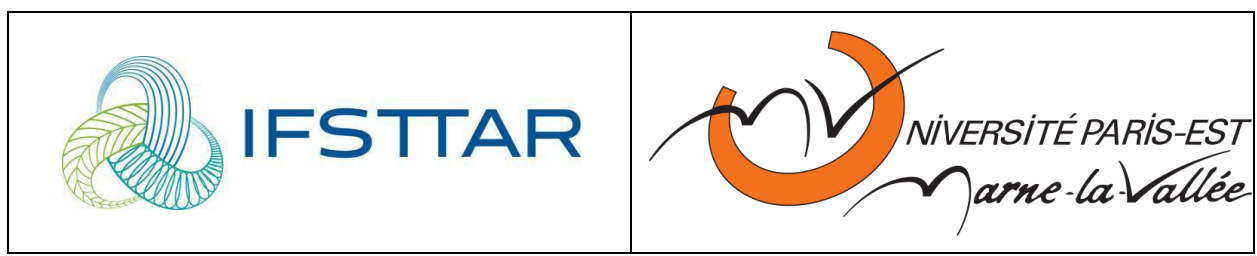

# Algebraic Approach for Performance Bound Calculus on Transportation Networks

(Road Network Calculus)

Nadir Farhi[*], Habib Haj-Salem[†] & Jean-Patrick Lebacque[‡]

Université Paris-Est, IFSTTAR, GRETTIA, F-93166 Noisy-le-Grand, France.

```
*nadir.farhi@ifsttar.fr (Corresponding author).
†habib.haj-salem@ifsttar.fr
‡jean-patrick.lebacque@ifsttar.fr
```

**Abstract**

We propose in this article an adaptation of the basic techniques of the deterministic *network calculus* theory to the road traffic flow theory. Network calculus is a theory based on min-plus algebra. It uses algebraic techniques to compute performance bounds in communication networks, such as maximum end-to-end delays and backlogs. The objective of this article is to investigate the application of such techniques for determining performance bounds on road networks, such as maximum bounds on travel times. The main difficulty to apply the network calculus theory on road networks is the modeling of interaction of cars inside one road, or more precisely the congestion phase. We propose a traffic model for a single-lane road without passing, which is compatible with the network calculus theory. The model permits to derive a maximum bound of the travel time of cars through the road. Then, basing on that model, we explain how to extend the approach to model intersections and large-scale networks.

**Keywords:** traffic flow theory, maximum travel time, network calculus, min-plus algebra.

## Notations

| | |
|---|---|
| $u(t)$ | inflow at time $t$. |
| $U(t)$ | cumulated inflow from time zero to time $t$. |
| $y(t)$ | outflow at time $t$. |
| $Y(t)$ | cumulated outflow from time zero to time $t$. |
| $\alpha$ | maximum arrival curve (time function). |
| $\beta$ | minimum service curve (time function). |
| $F$ | the set of time functions $\{f$ is non $-$ decreasing and $f(t) = 0, \forall\, t < 0\}$. |
| $\oplus$ | element-wise operation (min-plus addition in $F$).    $(f \oplus g)(t) := \min(f(t), g(t))$. |
| $*$ | min-plus convolution in $F$.    $(f * g)(t) := \inf_{0 \leq s \leq t}(f(s) + g(t - s))$. |
| $\oslash$ | min-plus de-convolution in $F$.    $(f \oslash g)(t) := \sup_{s \geq 0}(f(t + s) - f(s))$. |
| $\varepsilon$ | the zero element of the dioid $(F, \oplus, *)$.    $\varepsilon(t) = +\infty, \forall t \geq 0$. |
| $e$ | the unity element of the dioid $(F, \oplus, *)$.    $e(t) = +\infty, \forall t > 0$, and $e(0) = 0$. |
| $f^k$ | convolution power. $f^k = f * f * \ldots * f$    ($k$ times). |
| $B(t)$ | the backlog at time $t$.    $B(t) := U(t) - Y(t)$. |
| $d(t)$ | the virtual delay at time $t$.    $d(t) := \operatorname{Inf}\{h \geq 0, Y(t + h) \geq U(t)\}$. |
| $\gamma^p$ | a particular function in $F$.    $\gamma^p(t) = +\infty\; \forall t > 0$, and $\gamma^p(0) = p$. |
| $\delta^T$ | a particular function in $F$.    $\delta^T(t) = 0\; \forall t \leq T$, and $\delta^T(t) = +\infty\; \forall t > T$. |
| $[\text{expr}]^+$ | $\max(0, expr)$. |
| $q$ | car-flow. |
| $q_{max}$ | the maximum car-flow. |
| $q_i(t)$ | the car outflow from the $i^{\text{th}}$ section at time $t$. |
| $Q_i(t)$ | the cumulated car outflow from the $i^{\text{th}}$ section from time zero to time $t$. |
| $n_i$ | the number of cars in the $i^{\text{th}}$ section at time zero. |
| $n_{max}$ | the maximum number of cars that a section can contain. |
| $\rho$ | the average car-density in the road. |
| $\rho_c$ | the critical car-density. |
| $\rho_j$ | the jam car-density. |
| $v$ | the free car speed. |
| $w$ | the backward wave speed. |
| $\Delta x$ | the section length. |
| $q(\rho)$ | the car-flow function of the car-density (fundamental traffic diagram). |
| $\tau(\rho)$ | the travel time function of the car-density. |
| $\tau_{max}(\rho)$ | a maximum bound of the travel time, function of the car-density. |
| $\bigwedge_{k=1}^{m} f(k)$ | $\min_{1 \leq k \leq m} f(k)$. |

## 1  Introduction

The recent advances in information and communication technologies permit to obtain valuable information on the traffic state by means of probe vehicles (Amin & al, 2008), (Herrera & Bayen, 2008). The information is then either used to derive reliable traffic indicators, or sent (after required analyses, filtrations and reformulations) to connected vehicles via intelligent transportation equipments, in order to improve the traffic conditions. One of the most important traffic indicators that drivers need to receive in order to optimize their trip, is the travel time through the possible paths to their destinations. Even though the average value of the travel time estimation is determinant for drivers, its deviation may be very important in some cases. In order to evaluate the deviation of the travel time, one can determine either the probability distribution of the travel time, or minimum and maximum bounds for it. Several methods exist in the literature to estimate travel times (Coifman, 2002), (Claudel & Bayen, 2008), (Claudel, Hofleitner, Mignerey, & Bayen, 2008), (Ng, Szeto, & Waller, 2011).

We present in this article a traffic model that permits to derive a maximum bound of the travel time of cars passing through a single-lane road. The model is deterministic and uses algebraic techniques of the network calculus theory (a theory for performance bound calculus in communication and computer networks) (Chang, 2000), (Le Boudec & Thiran, 2001), (Jiang & Liu, 2008). The objective of our work is to adapt the algebraic approach of the *Network Calculus* theory to transportation networks.

We propose in this article a first step of applying the deterministic network calculus to determine the maximum bound of the travel time. The main contribution of this article is a traffic model that permits the derivation of that bound with a theoretic formula. Moreover, the maximum bound of the travel time is derived as a function of the average car-density on the road, and is then compared to the formula that gives the average travel time, basing on the same modeling. In section 2 we give a short review in the network calculus theory and in the min-plus algebra (Baccelli, Cohen, Olsder, & Quadrat, 1992), in order to fix the notations and the language. The model is presented in section 3. It is inspired from the cell transmission model (Daganzo, 1994), (Daganzo, 1995), and written on a single-lane road. The cumulative flows are seen as time signals and the car-dynamics is written algebraically as a min-plus linear system (Baccelli, Cohen, Olsder, & Quadrat, 1992). The impulse response of that system is then interpreted in term of guaranteed minimum service on the road. A maximum travel time of cars through the road is then derived from that guaranteed minimum service, as done basically in deterministic network calculus (Chang, 2000), (Le Boudec & Thiran, 2001). Thus, a formula giving a maximum bound for the travel time as a function of the average car-density in the road is obtained. In section 4, we propose a procedure for extending our approach to intersections and large transportation networks.

## 2   Review in network calculus

The procedure here is to consider a single-lane road as a server and apply the network calculus theory to determine a maximum bound of the travel time of cars through the road. In this short section we recall some basic results of the network calculus theory. In order to fix notations, a time function $f(t)$ (written with a small letter) expresses a given flow at time $t$, while $F(t)$ (written with a capital letter) denotes the cumulative flow $\int_0^t f(s)ds$. For an arrival flow $U$ arriving to a server that is empty of data at time zero, a maximum arrival curve $\alpha$ is associated in order to upper-bound the arrivals.

**Arrival curve.** $\alpha$ is a (maximum) arrival curve for $U$ if
$$U(t) - U(s) \leq \alpha(t-s), \forall\, 0 \leq s \leq t.$$

In the other side, a minimum service curve $\beta$ is associated to the server, in order to lower-bound the service. If the output from the server is denoted $Y$, then $\beta$ is defined as follows.

**Service curve.** $\beta$ is a (minimum) service curve for the server if
$$Y(t) \geq \min_{0 \leq s \leq t} \left(U(s) + \beta(t-s)\right), \forall t \geq 0.$$

Two indicators of the service performance are considered.
- The backlog $B(t)$ of data in the server at a given time $t$ is defined by
$$B(t) = U(t) - Y(t).$$

- The virtual delay $d(t)$ caused by the server at time $t$ is defined by
$$d(t) = \mathrm{Inf}\{h \geq 0, Y(t+h) \geq U(t)\}.$$

In the case where initial data $n_0$ is assumed in the server at time zero, the definition of the virtual delay remains correct by replacing the signal $Y$ by $Y - n_0$.

It is easy to see that arrival and service curves are not unique. In order to obtain good bounds on the backlog and on the virtual delay on a given server, it is necessary to consider "good" arrival and service curves. A good arrival (resp. service) curve is simply the minimal (resp. maximal) one; see (Le Boudec & Thiran, 2001) and (Chang, 2000) for more details.

We recall below a basic result of the deterministic network calculus (Le Boudec & Thiran, 2001) (Chang, 2000), that gives three bounds for a unique server. If $\alpha$ is an arrival curve for an arrival flow $U$ to a given server that offers a minimum service curve $\beta$, then we have

- The virtual delay is bounded as follows.

$$d(t) \leq \sup_{t \geq 0}\{\text{Inf}\{h \geq 0, \alpha(t) \leq \beta(t+h)\}\}, \quad \forall\, t \geq 0.$$

- The backlog is bounded as follows.

$$B(t) \leq \sup_{s \geq 0}(\alpha(s) - \beta(s)), \quad \forall t \geq 0.$$

- The curve $t \mapsto \sup_{s \geq 0}(\alpha(t+s) - \beta(s))$ is an arrival curve for the departure flow $Y$ from the server.

The maximum backlog and delay on a server are then given simply as the maximal vertical and horizontal distances between the arrival and the service curves; see Figure 1.

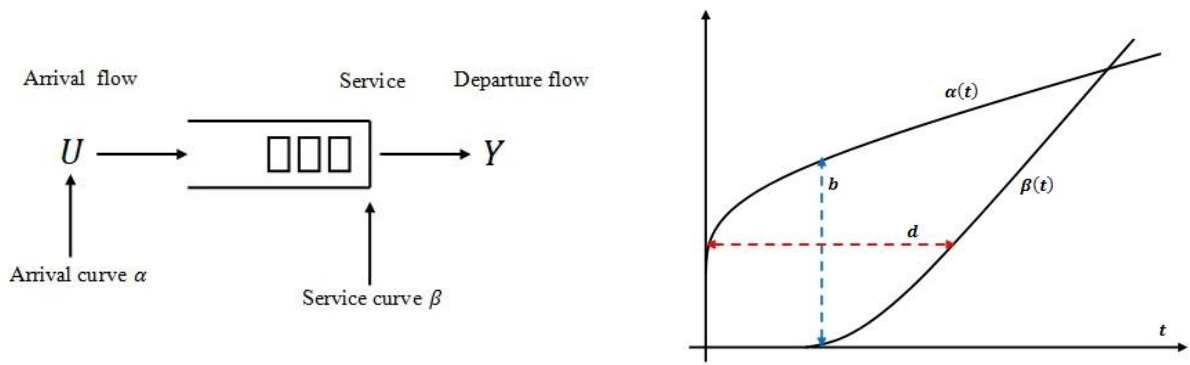

Figure 1. On the left side: Schema of the server. On the right side: The maximum delay $d$ and the maximum backlog $b$ determined graphically as the maximum horizontal and vertical distances between the arrival and the service curves, respectively.

## 3 Guaranteed service on a single-lane road (the model)

We present in this section an elementary model for the calculus of minimum guaranteed service for a single-lane road seen as a server. The objective of this modeling is to derive the maximum travel time of cars passing through the road, from the guaranteed minimum service of the road. The model is written on the cumulative car-flow variables $Q$ (see notations below). It is inspired from the cell transmission model (Daganzo, 1994), (Daganzo, 1995); see also (Lebacque, 1996). The min-plus linear system theory (Baccelli, Cohen, Olsder, & Quadrat, 1992) is then used to derive the minimum guaranteed service.

Let us consider a single-lane road where cars move without passing. In order to be able to fix the car-density, and to simplify the model, we consider a ring road; see Figure 2. The road is divided into $m$ sections of length $\Delta x$. The maximum number of cars on one section is denoted by $n_{max}$. The jam density $\rho_j$ on the road is then given by

$$\rho_j = \frac{n_{max}}{\Delta x}. \tag{1}$$

We use the following notations:

- $U(t)$ : The cumulated inflow of cars to the road from time zero up to time $t$.
- $Y(t)$ : The cumulated outflow of cars from the road, from time zero up to time $t$.
- $q_i(t)$ : The car outflow from the $i^{th}$ section (the section between positions $i\,\Delta x$ and $(i+1)\Delta x$) at time $t$.
- $Q_i(t) = \int_0^t q_i(s)ds$ : The cumulated car outflow from the $i^{th}$ section up to time $t$.
- $n_i$ : The number of cars in the $i^{th}$ section (the section between positions $i\,\Delta x$ and $(i+1)\Delta x$) at time zero.

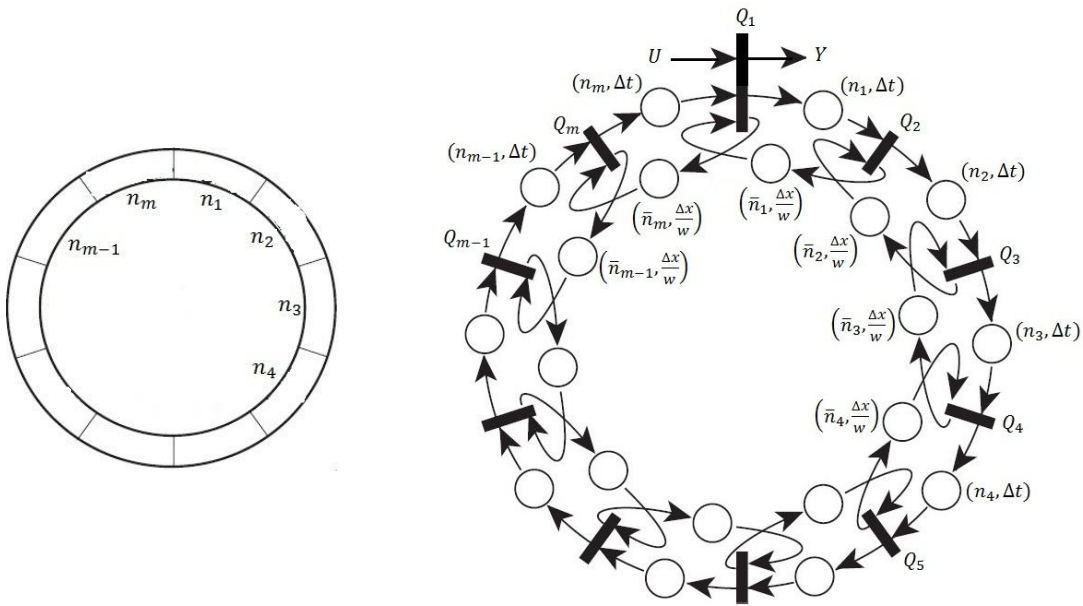

Figure 2. A single-lane ring road and the corresponding event-graph model.

Moreover, we assume the following fundamental traffic diagram for the road.

$$q(\rho) = \min\left(v\rho, w(\rho_j - \rho)\right), \tag{2}$$

Where $q, \rho, v, w$ denote respectively the car-flow, the car-density, the free car-speed and the backward wave speed on the road; see Figure 3. The maximum flow is then given by

$$q_{max} = \frac{\rho_j}{\frac{1}{v} + \frac{1}{w}}. \tag{3}$$

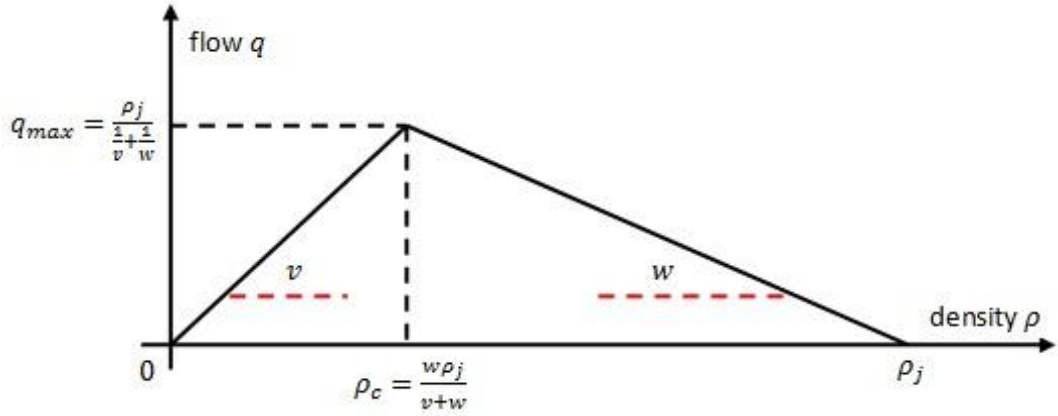

Figure 3. The fundamental traffic diagram.

Let us first recall the well known case (Farhi, 2008) of an autonomous single-lane ring road, with no entry $U$ and no exit $Y$. In this case the dynamics of the system is written as follows.

$$Q_1(t) = \min\{Q_m(t - \Delta x/v) + n_m, Q_2(t - \Delta x/w) + \bar{n}_1\}.$$
$$Q_i(t) = \min\{Q_{i-1}(t - \Delta x/v) + n_{i-1}, Q_{i+1}(t - \Delta x/w) + \bar{n}_i\}, \forall i \in \{2,3, \ldots, m-1\} \quad (4)$$
$$Q_m(t) = \min\{Q_{m-1}(t - \Delta x/v) + n_{m-1}, Q_1(t - \Delta x/w) + n_m\},$$

where $\bar{n}_k = \rho_j \Delta x - n_k = n_{max} - n_k$.

The system (4) is a discrete time event system for which the dynamics is well understood. Indeed, the system (4) can be represented by an event graph (a class of Petri nets), and its dynamics can be written linearly in min-plus algebra (Baccelli, Cohen, Olsder, & Quadrat, 1992). Since the event graph representing the system (4) is strongly connected, the asymptotic car flow on the road is then the same at every section, and is given by the minimum over the average weights of the graph circuits (Farhi, 2008). Three circuits are distinguished in the event graph of Figure 2.

- The interior circuit with the average growth rate
$$\frac{\sum_{i=1}^m n_i}{m \, \Delta x/v} = \frac{\sum_{i=1}^m n_i}{m \, \Delta x} v = v\rho.$$
- The exterior circuit with the average growth rate
$$\frac{\sum_{i=1}^m \bar{n}_i}{m \, \Delta x/w} = w(\rho_j - \rho).$$
- The loops over each section, with the average growth rate
$$\frac{n_{max}}{\Delta x/v + \Delta x/w} = \frac{\rho_j \Delta x}{\Delta x/v + \Delta x/w} = q_{max}.$$

The asymptotic car-flow is then given by

$$q(\rho) = q_i(\rho) = \min\{v\rho, -w(\rho - \rho_j)\}, \forall i.$$

and retrieves then the fundamental diagram (2) assumed in the model.

The average travel time through the road is then given by

$$\tau(\rho) = m\Delta x \frac{\rho}{q(\rho)} = \max\left\{m\frac{\Delta x}{v}, m\frac{\Delta x}{w}\frac{\rho}{(\rho_j - \rho)}\right\}.$$

In the following, we are interested on the maximum travel time through the road, rather than the average travel time. For that, we introduce an inflow and outflow of cars on the road. We assume that cars enter into the road from the entry of section 1, and exit it from the same point after passing through the road. The dynamics (4) are then rewritten as follows:

$$\begin{aligned}
Q_1(t) &= \min\{U(t), Q_m(t - \Delta x/v) + n_m, Q_2(t - \Delta x/w) + \bar{n}_1\}. \\
Q_i(t) &= \min\{Q_{i-1}(t - \Delta x/v) + n_{i-1}, Q_{i+1}(t - \Delta x/w) + \bar{n}_i\}, \forall\, i \in \{2,3,\ldots,m-1\} \\
Q_m(t) &= \min\{Q_{m-1}(t - \Delta x/v) + n_{m-1}, Q_1(t - \Delta x/w) + \bar{n}_m\}. \\
Y(t) &= \min\{Q_m(t - \Delta x/v) + n_m, Q_2(t - \Delta x/w) + \bar{n}_1\}.
\end{aligned} \quad (5)$$

Let us notice the difference between the dynamics of $Q_1$ which is conditioned by the inflow $U$, and the dynamics of $Y$ which is not conditioned by the inflow $U$.

We consider the same modeling approach as above, but here we will see the variables as signals; see (Baccelli, Cohen, Olsder, & Quadrat, 1992). We first show that the dynamic system (5) is linear in a certain algebraic structure, and then use this linearity to derive a guaranteed service of the road seen as a car server.

Let us consider the set of time functions
$$F = \{f \text{ is non} - \text{decreasing and } f(t) = 0, \forall\, t < 0\},$$
endowed with the following two operations.

- Addition (element-wise minimum): $(f \oplus g)(t) := \min\{f(t), g(t)\}, \forall t \geq 0$.
- Product (min-plus convolution): $(f * g)(t) := \min_{0 \leq s \leq t}\{f(s) + g(t-s)\}, \forall t \geq 0$.

We then obtain a dioid structure $(F, \oplus, *)$ (Baccelli, Cohen, Olsder, & Quadrat, 1992) (Le Boudec & Thiran, 2001). The zero element of that dioid is denoted by $\varepsilon$ and defined by $\varepsilon(t) = +\infty, \forall t \geq 0$, while the unity element is denoted by $e$ and defined by

$$e(t) = \begin{cases} 0 \text{ if } t = 0. \\ +\infty \text{ for } t \geq 0. \end{cases}$$

In addition, we consider the following notations:

- Power: $f^k$ denotes the product (min-plus convolution) of $f$ with itself $k$ times.
- Additive closure: $f^*(t) := e \oplus f \oplus f^2 \oplus f^3 \oplus \ldots$
- De-convolution: $(f \oslash g)(t) := \sup_{s \geq 0}(f(t+s) - g(s))$.

and the two particular signals $\gamma^p$ (the gain signal) and $\delta^T$ (the shift signal) in $F$ (see Appendix A for more details on those signals).

$$\gamma^p(t) = \begin{cases} p \text{ if } t = 0 \\ +\infty \text{ for } t > 0 \end{cases} \quad \text{and} \quad \delta^T(t) = \begin{cases} 0 \text{ if } t \leq T \\ +\infty \text{ otherwise} \end{cases}$$

We notice that by using the notation of the min-plus convolution, the definitions of arrival and service curves can be rewritten as follows.
- $\alpha$ is an arrival curve for $U$ if $U \leq \alpha * U$.
- $\beta$ is a (minimum) service curve for the server if $Y \geq \beta * U$.

Now, since the cumulative flows $U, Q_i, 1 \leq i \leq m$, and $Y$ are time functions (or signals) in $F$, then by using the notations of addition and product in $F$, the system (5) is linear according to those operations, and is written as follows[1].

$$
\begin{aligned}
Q_1 &= \gamma^{n_m} \delta^{\Delta x/v} Q_m \oplus \gamma^{\bar{n}_1} \delta^{\Delta x/w} Q_2 \oplus U. \\
Q_i &= \gamma^{n_{i-1}} \delta^{\Delta x/v} Q_{i-1} \oplus \gamma^{\bar{n}_i} \delta^{\Delta x/w} Q_{i+1}, \quad \forall i \in \{2,3,\ldots,m-1\}. \\
Q_m &= \gamma^{n_{m-1}} \delta^{\Delta x/v} Q_{m-1} \oplus \gamma^{\bar{n}_m} \delta^{\Delta x/w} Q_1. \\
Y &= \gamma^{n_m} \delta^{\Delta x/v} Q_m \oplus \gamma^{\bar{n}_1} \delta^{\Delta x/w} Q_2.
\end{aligned}
\quad (6)
$$

Moreover, since we are not only interested in the average quantities, but on the maximum bounds, we include the initial conditions of the system (5) (the signals are null at time zero). To include those conditions, it is sufficient to add (min-plus addition) to each signal of system (6) the signal unity $e$. Then, the dynamics (6) are written as follows.

$$
\begin{aligned}
Q_1 &= \gamma^{n_m} \delta^{\Delta x/v} Q_m \oplus \gamma^{\bar{n}_1} \delta^{\Delta x/w} Q_2 \oplus U \oplus e. \\
Q_i &= \gamma^{n_{i-1}} \delta^{\Delta x/v} Q_{i-1} \oplus \gamma^{\bar{n}_i} \delta^{\Delta x/w} Q_{i+1} \oplus e, \quad \forall i \in \{2,3,\ldots,m-1\}. \\
Q_m &= \gamma^{n_{m-1}} \delta^{\Delta x/v} Q_{m-1} \oplus \gamma^{\bar{n}_m} \delta^{\Delta x/w} Q_1 \oplus e. \\
Y &= \gamma^{n_m} \delta^{\Delta x/v} Q_m \oplus \gamma^{\bar{n}_1} \delta^{\Delta x/w} Q_2 \oplus e.
\end{aligned}
\quad (7)
$$

The system (7) is then written as follows.

$$
\begin{aligned}
Q &= A * Q \oplus B * U \oplus E \\
Y &= C * Q \oplus e
\end{aligned}
\quad (8)
$$

where $Q = (Q_1 \ Q_2 \ \cdots \ Q_m)'$,

$$
A = \begin{pmatrix}
\varepsilon & \gamma^{\bar{n}_1}\delta^{\Delta x/w} & \varepsilon & \cdots & \varepsilon & \gamma^{n_m}\delta^{\Delta t} \\
\gamma^{n_1}\delta^{\Delta t} & \varepsilon & \gamma^{\bar{n}_2}\delta^{\Delta x/w} & \varepsilon & \cdots & \varepsilon \\
\varepsilon & \gamma^{n_2}\delta^{\Delta t} & \varepsilon & \gamma^{\bar{n}_3}\delta^{\Delta x/w} & \varepsilon & \vdots \\
\vdots & \ddots & \ddots & \ddots & \ddots & \varepsilon \\
\varepsilon & \cdots & \varepsilon & \gamma^{n_{m-2}}\delta^{\Delta t} & \varepsilon & \gamma^{\bar{n}_{m-1}}\delta^{\Delta x/w} \\
\gamma^{\bar{n}_m}\delta^{\Delta x/w} & \varepsilon & \cdots & \varepsilon & \gamma^{n_{m-1}}\delta^{\Delta t} & \varepsilon
\end{pmatrix}, \quad
B = \begin{pmatrix} e \\ \varepsilon \\ \vdots \\ \varepsilon \\ \varepsilon \end{pmatrix},
$$

$C = (\varepsilon \ \gamma^{\bar{n}_1}\delta^{\Delta x/w} \ \varepsilon \ \cdots \ \varepsilon \ \gamma^{n_m}\delta^{\Delta t})$ and $E = (e \ e \ \cdots \ e)'$.

The system (8) is a min-plus linear system. A basic result of the min-plus system theory (Baccelli, Cohen, Olsder, & Quadrat, 1992) gives then the impulse response of system (8):

---

[1] Note that, as in the standard algebra, the product operation is sometimes just not symbolized (that is to say that $f * g$ can simply be written $fg$.)

$$Y = CA^*(BU \oplus E) \oplus e = CA^*BU \oplus CA^*E \oplus e. \tag{9}$$

Note that an impulse response $Y = CA^*BU$ would be obtained if system (6) is considered instead of system (7) (that is if the initial conditions of system (6) are not taken into account). In that case, we would conclude directly that $CA^*B$ is a minimum service curve of the road seen as a server (since we have $Y \geq (CA^*B) * U$), and derive a maximum bound for the travel time through the road. But here, as mentioned above, it is necessary to take into account the initial conditions of the dynamical system, since we are interested in the maximum bound of the travel time through the road, rather than the average travel time. In order to be able to derive maximum bounds from the formula (9), as done from a service curve, we propose the following extension.

**Minimum service couple**

We consider here the case where the curve service is given with an additional affine term. More precisely, we say that $(\beta, \lambda)$ is a *service couple* for a server if

$$Y \geq \beta * U \oplus \lambda.$$

Then we can easily check (see Appendix B), that to obtain the three bounds given above, for that case, it is sufficient to replace $\beta$ with $\beta \oplus \lambda$. That is to say that:

- The maximum backlog is bounded as follows.
$$B(t) \leq \sup_{s \geq 0}(\alpha(s) - (\beta \oplus \lambda)(s)), \quad \forall t \geq 0.$$
- The maximum delay is bounded as follows.
$$d(t) \leq \sup_{t \geq 0}\{\text{Inf}\{h \geq 0, \alpha(t) \leq (\beta \oplus \lambda)(t + h)\}\}, \quad \forall t \geq 0.$$
- An arrival curve for the departure flow is $\alpha \oslash (\beta \oplus \lambda)$. That is
$$Y \leq (\alpha \oslash (\beta \oplus \lambda)) * Y.$$

Note that the curve $(\beta \oplus \lambda)$ is not necessarily a minimu service curve for the server.

The two first bounds are then given by the maximum vertical and horizontal distances between the curves $\alpha$ and $\beta \oplus \lambda$.

**Theorem 1.** A minimum service couple for the single-lane road, seen as a server, is $(\beta, \lambda)$ given by

$$\beta = [\gamma^{-m\rho\Delta x} * a^*]^+$$

$$\lambda = \beta \oplus \left[\gamma^{-m\rho\Delta x}\left(\bigwedge_{k=1}^{m-1} \gamma^{[m\rho-k\rho_j]^+\Delta x}\delta^{(m-k)\Delta x/v} \oplus \bigwedge_{k=1}^{m-1} \gamma^{[k\rho_j-m\rho]^+\Delta x}\delta^{k\Delta x/w} \oplus e\right)\right]^+$$

where

$$a = \gamma^{m\rho\Delta x}\delta^{m\Delta x/v} \oplus \gamma^{\rho_j\Delta x}\delta^{\Delta x/v+\Delta x/w} \oplus \gamma^{m(\rho_j-\rho)\Delta x}\delta^{m\Delta x/w}$$

**Proof.** The system (8) is an affine min-plus system. We then have $Y = CA^*BU \oplus CA^*E \oplus e$. We need to compute $A^*$. Several methods can be used to compute $A^*$. By definition of $A^*$, one can simply compute $A^2, A^3, \ldots$etc, then deduce $A^* = \bigwedge_{k \geq 0} A^k$ by simplifying all the terms. We can easily check that $A^*$ is given as follows.

$$(A^*)_{ij} = \begin{cases} a^* \text{ if } i = j \\ a^* \left( \gamma^{\sum_{k=j}^{i-1} n_k} \delta^{(i-j)\Delta t} \oplus \gamma^{\sum_{k=i}^{j-1} \bar{n}_k} \delta^{(j-i)\Delta x/w} \right) \text{ if } i \neq j, \text{with k cyclic in } \{1,2,\ldots,m\} \end{cases}.$$

In order to compute $CA^*B$, we need only the first column of $A^*$, since only the first entry of $B$ is not null. Thus, $A^*B = (A^*)_{\cdot 1}$ and is given by

$$A^*B = a^* * \left( e \quad (\gamma^{n_1} \delta^{\Delta t} \oplus \gamma^{\sum_{k=2}^{m} \bar{n}_k} \delta^{(m-1)\Delta x/w}) \quad \cdots \quad (\gamma^{\sum_{k=1}^{m-1} n_k} \delta^{(m-1)\Delta t} \oplus \gamma^{\bar{n}_m} \delta^{\Delta x/w}) \right).$$

Therefore, $CA^*B = a^* * a$.

Similarly, we can check that

$$CA^*E = a^* * a \oplus \bigwedge_{k=1}^{m-1} \gamma^{\sum_{i=k+1}^{m} n_i} \delta^{(m-k)\Delta x/v} \oplus \bigwedge_{k=1}^{m-1} \gamma^{\sum_{i=1}^{k} \bar{n}_i} \delta^{k\,\Delta x/w}$$

Then, it is not difficult to check that

$$\sum_{i=k+1}^{m} n_i \geq [m\rho - k\rho_j]^+ \Delta x.$$
$$\sum_{i=1}^{k} \bar{n}_i \geq [k\rho_j - m\rho]^+ \Delta x.$$

Now, since we are interested in the car outflow from the road that comes from the car inflow, without counting the cars being in the road at time zero (and stay there all time because the road is circular), we need to express the variable $Z = [Y - \sum_{i=1}^{m} n_i]^+ = [\gamma^{-m\rho \Delta x} Y]^+$ in function of $U$. We then conclude that the couple of curves $(\gamma^{-m\rho \Delta x} CA^*B, \gamma^{-m\rho \Delta x} (CA^*E \oplus e))$ is a minimum service couple for the ring road (see property (P4) in Appendix A).

Finally, we have $[\gamma^{-m\rho \Delta x} * a * a^*]^+ = [\gamma^{-m\rho \Delta x} * a^*]^+$, since $(a * a^*)(t) = a^*(t), \forall t > 0$ and $(a * a^*)(0) = a(0) < m\rho \Delta x$ (see properties (P2) and (P5) in Appendix A).
∎

In order to show the shape of the minimum service couple given in Theorem 1, let us take an academic example.

**Example 1.** Let $m = 6, \Delta x = 1, v = 1, w = 1/2, \rho_j = 1$. Then we have $n_{max} = 1, \rho_c = 1/3$ and $q_{max} = 1/3$. In addition, we take three cases, where we vary the average car-density on the road. Let us first notice that $a^*$ can also be written as follows (see property (P3) in Appendix A).

$$a^* = (\gamma^{m\rho \Delta x} \delta^{m\Delta x/v})^* * (\gamma^{\rho_j \Delta x} \delta^{\Delta x/v + \Delta x/w})^* * (\gamma^{m(\rho_j - \rho)\Delta x} \delta^{m\Delta x/w})^*$$

Then we have the three cases:

- $\sum_{i=1}^{m} n_i = 1$, that is $\rho = 1/6 < \rho_c$. Then $a = \gamma^1 \delta^6 \oplus \gamma^5 \delta^{12}$ and $a^* = (\gamma^1 \delta^6)^*$.

Hence the minimum service couple $(\beta, \lambda)$ is given by
$$\beta = [\gamma^{-1}(\gamma^1 \delta^6)^*]^+ = \gamma^{-1} * (\gamma^1 \delta^6)^*.$$
$$\lambda = \beta \oplus \gamma^0 \delta^5 \oplus \gamma^1 \delta^6 \oplus \gamma^2 \delta^8 \oplus \gamma^3 \delta^{10} = \beta.$$

- $\sum_{i=1}^{m} n_i = 2$, that is $\rho = 1/3 = \rho_c$. Then $a = \gamma^1 \delta^3$ and $a^* = (\gamma^1 \delta^3)^*$. Hence the minimum service couple $(\beta, \lambda)$ is given by
$$\beta = [\gamma^{-2}(\gamma^1 \delta^3)^*]^+.$$
$$\lambda = \beta \oplus \gamma^0 \delta^8 \oplus \gamma^1 \delta^{10}.$$

- $\sum_{i=1}^{m} n_i = 3$, that is $\rho = 1/2 > \rho_c$. Then $a = \gamma^1 \delta^3 \oplus \gamma^3 \delta^{12}$, $a^* = (\gamma^1 \delta^3)^*(\gamma^3 \delta^{12})^*$.

  Hence the minimum service couple $(\beta, \lambda)$ is given by
$$\beta = [\gamma^{-3}(\gamma^1 \delta^3)^*(\gamma^3 \delta^{12})^*]^+.$$
$$\lambda = \beta \oplus \gamma^0 \delta^{10} = \beta.$$

The minimum service couples (given in Theorem 1) corresponding to each of the three cases, are shown in Figure 4.

∎

We give below a corollary of Theorem 1, where by relaxing the minimum service couple given in Theorem 1, we obtain practical formulas.

**Corollary 1.** A minimum service couple $(\beta, \lambda)$ for the road is given by
$$\beta(t) = q(\rho)[t - \tau(\rho)]^+$$
$$\lambda(t) = \beta(t) \oplus w\rho_j \left[t - \frac{2m\rho}{\rho_j} \frac{\Delta x}{w}\right]^+$$

where $q(\rho)$ and $\tau(\rho)$ are the average flow and average travel time given by the fundamental diagram (given in section 2).

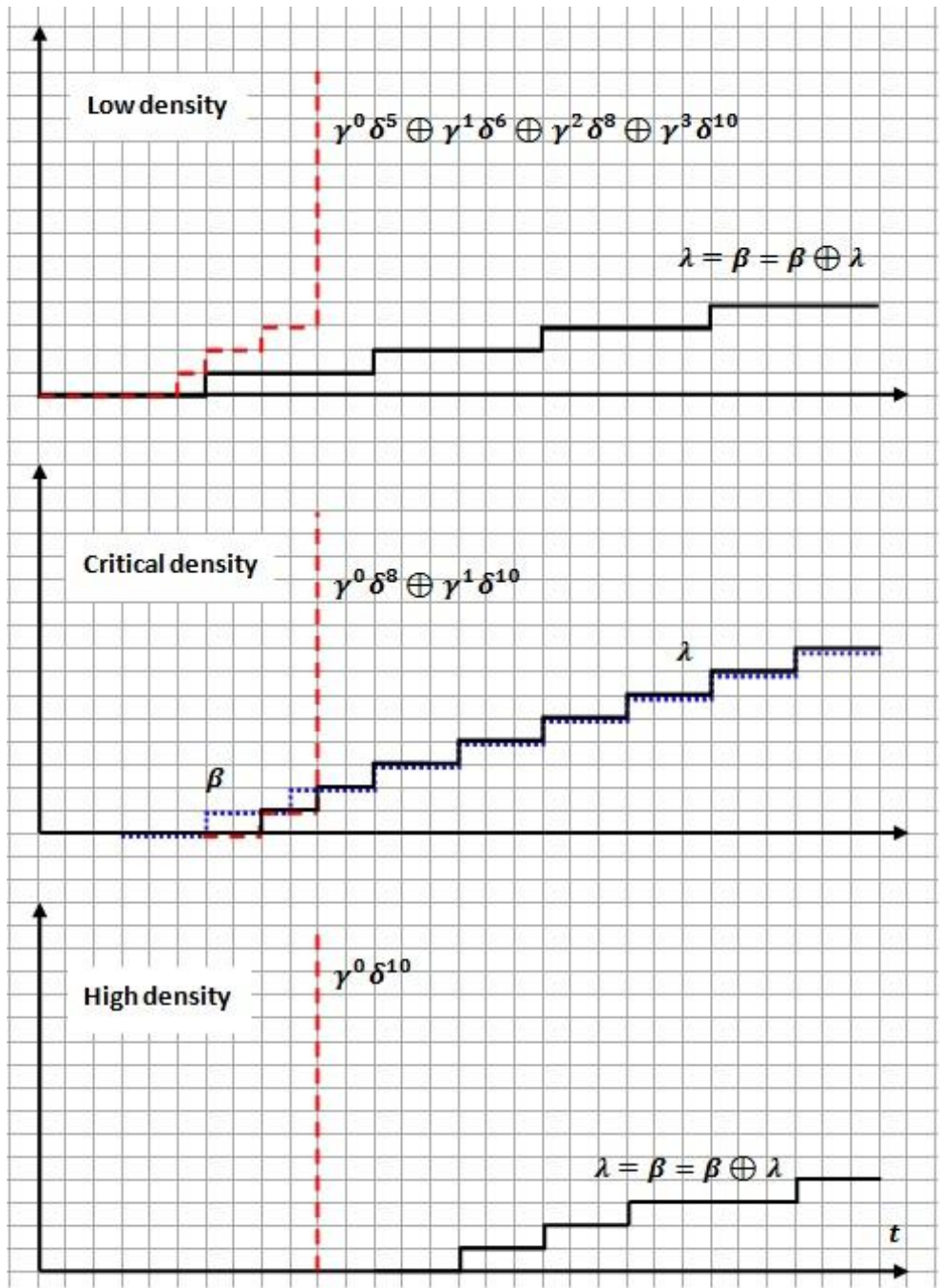

Figure 4. Calculus of the service couple for the single-lane ring road in three cases of low, critical and high car-density.

**Proof.**

1. From the curve $\beta$ given by Theorem 1, we have

   - $\left(\gamma^{\sum_{i=1}^{m} n_i} \delta^{m\frac{\Delta x}{v}}\right)^* \geq \rho v \, t,$
   - $\left(\gamma^{n_{max}} \delta^{\frac{\Delta x}{v}+\frac{\Delta x}{w}}\right)^* \geq q_{max} \, t,$
   - $\left(\gamma^{\sum_{i=1}^{m} \bar{n}_i} \delta^{m\frac{\Delta x}{w}}\right)^* \geq (\rho_j - \rho) w t.$

Then $\left(\gamma^{\sum_{i=1}^{m} n_i} \delta^{m\frac{\Delta x}{v}}\right)^* * \left(\gamma^{n_{max}} \delta^{\frac{\Delta x}{v}+\frac{\Delta x}{w}}\right)^* * \left(\gamma^{\sum_{i=1}^{m} \bar{n}_i} \delta^{m\frac{\Delta x}{w}}\right)^* \geq \min\{\rho v, q_{max}, (\rho_j - \rho)w\} \, t.$

But since $\forall \rho, \min\{\rho v, (\rho_j - \rho)w\} \geq q_{max}$, then

$$\left(\gamma^{\sum_{i=1}^{m} n_i} \delta^{m\frac{\Delta x}{v}}\right)^* * \left(\gamma^{n_{max}} \delta^{\frac{\Delta x}{v}+\frac{\Delta x}{w}}\right)^* * \left(\gamma^{\sum_{i=1}^{m} \bar{n}_i} \delta^{m\frac{\Delta x}{w}}\right)^* \geq \min\{\rho v, (\rho_j - \rho)w\} \, t.$$

Hence

$$\left[\left(\gamma^{\sum_{i=1}^{m} n_i} \delta^{m\frac{\Delta x}{v}}\right)^* * \left(\gamma^{n_{max}} \delta^{\frac{\Delta x}{v}+\frac{\Delta x}{w}}\right)^* * \left(\gamma^{\sum_{i=1}^{m} \bar{n}_i} \delta^{m\frac{\Delta x}{w}}\right)^* - \rho m \Delta x\right]^+$$
$$\geq \min\{\rho v, (\rho_j - \rho)w\}\left(t - \frac{\rho m \Delta x}{\min\{\rho v, (\rho_j - \rho)w\}}\right)^+$$
$$\geq \min\{\rho v, (\rho_j - \rho)w\}\left(t - \max\left\{m\frac{\Delta x}{v}, \frac{\rho}{\rho_j - \rho} m \frac{\Delta x}{w}\right\}\right)^+.$$

2. From the curve $\lambda$ given in theorem 1, we have

$$\bigwedge_{k=1}^{m-1} \gamma^{[m\rho - k\rho_j]^+ \Delta x} \delta^{(m-k)\Delta x / v}$$
$$\geq \max\left\{v\rho_j \left[t - \left(m - \frac{m\rho}{\rho_j}\right)\frac{\Delta x}{v}\right]^+, \gamma^{[m\rho - \rho_j]^+ \Delta x} \delta^{(m-1)\Delta x / v}\right\}$$

then

$$\gamma^{-m\rho \Delta x} \bigwedge_{k=1}^{m-1} \gamma^{[m\rho - k\rho_j]^+ \Delta x} \delta^{(m-k)\Delta x / v} \geq v\rho_j \left[t - m\frac{\Delta x}{v}\right]^+.$$

In the other hand, we have

$$\bigwedge_{k=1}^{m-1} \gamma^{[k\rho_j - m\rho]^+ \Delta x} \delta^{k \Delta x / w} \geq \max\left\{w\rho_j \left[t - \frac{m\rho}{\rho_j}\frac{\Delta x}{w}\right]^+, \gamma^{[(m-1)\rho_j - m\rho]^+ \Delta x} \delta^{(m-1)\Delta x / w}\right\}$$

then

$$\gamma^{-m\rho \Delta x} \bigwedge_{k=1}^{m-1} \gamma^{[k\rho_j - m\rho]^+ \Delta x} \delta^{k \Delta x / w} \geq w\rho_j \left[t - \frac{2m\rho}{\rho_j}\frac{\Delta x}{w}\right]^+.$$

3. Finally, it is easy to check that

$$\beta(t) \leq v\rho_j \left[t - m\frac{\Delta x}{v}\right]^+$$

since $q(\rho) \leq v\rho_j$ and $\tau(\rho) \geq m \Delta x / v$. ∎

Let us notice that the term $w\rho_j \left[t - \frac{2m\rho}{\rho_j}\frac{\Delta x}{w}\right]^+$ is more important than the term $\beta(t)$ in the service couple given in Corollary 1. Indeed the curve $\beta(t)$ gives simply the average service of the road, since $q(\rho)$ is the average car-flow and $\tau(\rho)$ is the average travel time on the road.

**Corollary 2.** If $\alpha(t) = \sigma + rt$ is an arrival curve for the inflow $U$ to the road, then the following bounds are guaranteed.

- A maximum bound for the travel time of cars through the road is
$$\tau_{max} = \max\left\{\tau + \frac{\sigma}{q}, \frac{2m\rho}{\rho_j}\frac{\Delta x}{w} + \frac{\sigma}{w\rho_j}\right\}.$$
- A maximum bound for the number of cars waiting at the entry of the road (not confuse with the number of cars queuing at the exit of the road) is
$$b_{max} = \max\left\{\sigma + r\tau, \sigma + r\frac{2m\rho}{\rho_j}\frac{\Delta x}{w}\right\}.$$
- An arrival curve for the departure flow from the road is
$$\bar{\alpha}(t) = b_{max} + rt.$$

**Proof.** It is well known (Le Boudec & Thiran, 2001) that if a flow with an arrival curve $\alpha(t) = \sigma + rt$ is served in a server with a minimum service curve $\beta(t) = R(t - T)^+$, then the maximum delay is $T + \sigma/R$, the maximum backlog is $\sigma + rT$, and the curve $\sigma + rT + rt$ is an arrival curve for the departure flow. The result is then obtained by adapting these bounds to the case of couple service instead of minimum service curve (see Appendix B), and by using the couple service given in Corollary 1.

∎

In the road traffic, it is probably not interesting to assume arrival curves with no null $\sigma$. One may simply estimate the arrival flow rate (linear arrival curve) and the car-density on the road, at a given time, and want to determine the maximum three bounds at the considered time instant. In this case, we simply have from Corollary 2:

$$\tau_{max}(\rho) = \max\left\{\tau(\rho), \frac{2m\rho}{\rho_j}\frac{\Delta x}{w}\right\}$$
$$= \max\left\{\frac{1}{v}, \frac{\rho}{\rho_j - \rho}\frac{1}{w}, \frac{2\rho}{\rho_j}\frac{1}{w}\right\} m\Delta x, \tag{10}$$

and $b_{max} = r\,\tau_{max}$, and $\bar{\alpha}(t) = r(t + \tau_{max})$.

The formula (10) tells that the maximum travel time through the road is greater than the average travel time only in the car-density interval $[(\rho_j/2)(w/v), \rho_j/2]$ when $w < v$. That is to say that

$$\frac{\rho_j}{2}\frac{w}{v} < \rho < \frac{\rho_j}{2} \Leftrightarrow \tau_{max}(\rho) > \tau(\rho). \tag{11}$$

In Figure 5, we show the average and the maximum travel times in function of the car-density in the case where $w < v$.

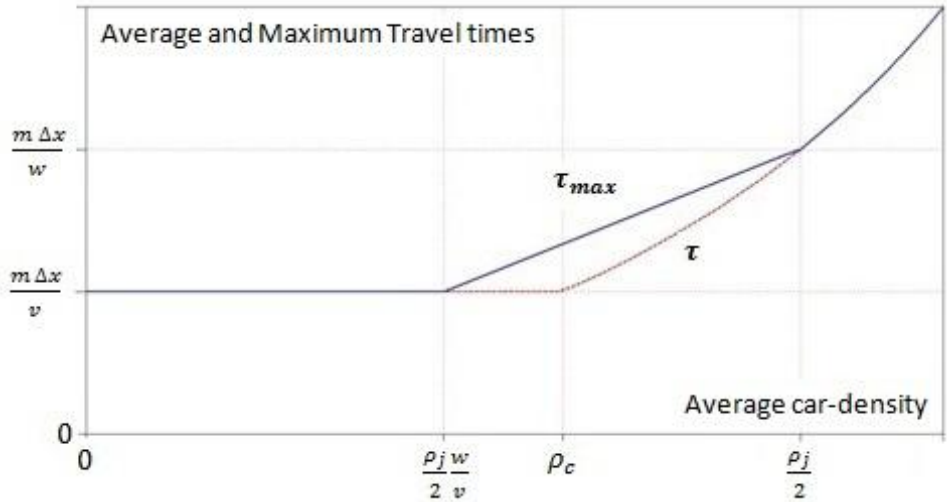

Figure 5. The average and the maximum travel times through the road in the case where $w = v/2$.

Let us notice that the choice of the traffic flow model (cell transmission model) and that of the fundamental traffic diagram (triangular one) are very important. Indeed these two choices permitted to write the car-dynamics on the road linearly in the min-plus algebra. This is necessary because the service couple is derived as an impulse response of a min-plus linear system.

## 4  Road network calculus

We give in this section the first ideas on how to use the results obtained on a single-lane road to extend the approach to intersections and whole transportation networks. First we notice that for multi-lane roads, one may just use one fundamental diagram for all lanes and apply the same approach as done for a single-lane road. We assume here that the same model and results are used for multi-lane roads. However, other models similar to the one given above can also be developed for multi-lane traffic.

Before presenting the procedure of extending our approach to intersections and networks, we need to recall and extend an important result of the deterministic network calculus on the series composition of servers (it shows in particular the power of the algebraic approach against other approaches). The result tells that a minimum service curve of the series composition of two servers guaranteeing $\beta_1$ and $\beta_2$ as minimum service curve for each of them, is simply the curve $\beta_1 * \beta_2$. The result can easily be proved by using the associative property of the min-plus convolution; see (Chang, 2000) (Le Boudec & Thiran, 2001) for more details. Moreover, since the min-plus convolution is also commutative, then the order of the composition of the two servers is not important.

For our model, we need to have a similar result for a composition of servers offering minimum service couples (rather than minimum service curves). It is easy to see that the composition of two servers offering two service couples $(\beta_1, \lambda_1)$ and $(\beta_2, \lambda_2)$ guarantees a service couple $(\beta_1 * \beta_2, \beta_1 * \lambda_2 \oplus \lambda_1)$. Indeed from $Y = \beta_1 * Z \oplus \lambda_1$ and $Z = \beta_2 * U \oplus \lambda_2$, we get $Y = (\beta_1 * \beta_2) * U \oplus (\beta_1 * \lambda_2 \oplus \lambda_1)$.

In communication and computer networks, one determines first the (residual) guaranteed service for all (input, output) couples through every switching router. Then, it is sufficient to compose (with a min-plus convolution) all the guaranteed services through all the arcs of a given path, to determine the service guaranteed through the whole path. Finally, one determines the maximum end-to-end delays on a communication network by calculating the maximum delay on each (origin-destination) path on the network, by simply using its guaranteed service curve, and considering the whole path as an elementary server. The residual guaranteed service calculus on a given input-output couple of a given router takes into account the control policy set in that router; see (Chang, 2000) and (Le Boudec & Thiran, 2001) for more details.

We think that we can proceed similarly for transportation networks. Indeed, the main difference between data traffic in communication networks and car-traffic in transportation networks is the interaction between particles (drivers observe reaction times contrary to data packets), expressed in the fundamental diagram of the road. Since this difference is already taken into account in the one road model (presented above), then to extend the approach to complicated transportation networks, it remains the adaptation of the residual guaranteed service calculus on input-output couples of routers to apply in intersections, for transportation control policies. Then one only needs to compose elementary road services and residual guaranteed services, to obtain guaranteed services on whole paths. Maximum travel times can then be derived similarly. We notice here that one of the most difficult issues to solve, in order to extend the approach presented in this article to big networks, is the presence of cyclic dependencies of inflows arriving to one intersection, even though some elementary results exist to deal with that issue (Chang, 2000).

Let us explain the calculus of maximum bounds of travel times in a tree-like network. Let us consider the transportation network of Figure 6.

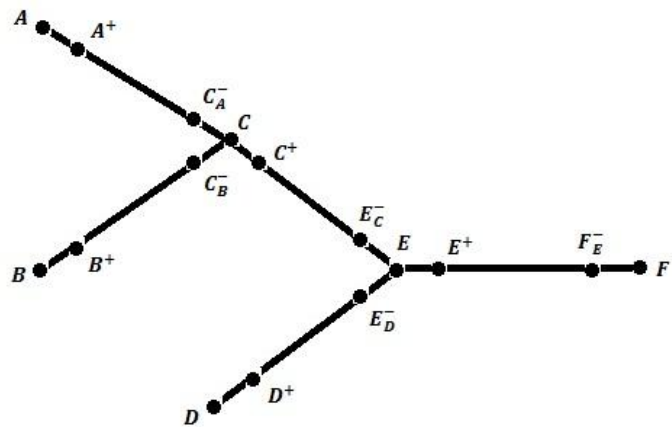

Figure 6. Tree-like network. The car-traffic goes from the left side to the right side.

In Figure 6, the notation $A^+$ indicates some point downstream of intersection $A$. The notation $C_A^-$ indicates some point upstream of intersection $C$ in the direction of intersection $A$. All other notations are interpreted similarly. We assume that we have the fundamental diagrams of the roads $(A^+, C_A^-)$, $(B^+, C_B^-)$, $(C^+, E_C^-)$, $(D^+, E_D^-)$ and $(E^+, F_E^-)$. The approach for calculating a maximum bound for the travel time from $A$ to $F$ in Figure 6 is the following.

1. Determine the guaranteed service on each of the roads $(A^+, C_A^-)$, $(B^+, C_B^-)$, $(C^+, E_C^-)$, $(D^+, E_D^-)$ and $(E^+, F_E^-)$, independent of the intersections. We simply use the fundamental diagrams of the roads, and the model presented above.
2. Calculate the residual services on the merge $C$:
   a. The service guaranteed on $C$ for the aggregate inflows coming from $C_A^-$ and $C_B^-$ is assumed to be simply the guaranteed service on the road $(C^+, E_C^-)$.
   b. The control policy on $C$ is assumed to be known.
   c. From these two information (a and b), calculate the residual guaranteed services for the flows $(C_A^-, C^+)$ and $(C_B^-, C^+)$. This shall be done by adapting well known results of network calculus theory on the calculus of residual guaranteed service. (This is not yet done).
3. Calculate the residual services on the merge $E$. By the same method as in 2. Calculate the residual guaranteed services for the flows $(E_C^-, E^+)$ and $(E_D^-, E^+)$.
4. Finally, the guaranteed service for the flow $(A^+, F_E^-)$, or simply $(A, F)$, is given by the series composition of the services $(A^+, C_A^-)$, $(C_A^-, C^+)$, $(C^+, E_C^-)$, $(E_C^-, E^+)$ and $(E^+, F_E^-)$, respectively. The guaranteed services for $(B, F)$ and $(D, F)$ are obtained similarly.

The maximum travel times from $A$ (resp. $B$ and $D$) to $F$ are then derived from the guaranteed services $(A, F)$, $(B, F)$ and $(D, F)$ respectively, as done on a single road.

## Conclusion

We presented in this article a network calculus traffic model that permits the derivation of a maximum bound for the travel time of cars passing through a single-lane road. An important advantage of this model is its algebraic formulation which is very powerful comparing to other formulations (e.g. the series composition). Even though the model is basic and elementary, its developments and extensions may be promising. Our future work shall be on the realization of the model extension process presented in section 5. We shall also demonstrate the effectiveness of our approach by performing numerical investigations on effective data sets.

## Appendix A (Details on the signals $\gamma^p$ and $\delta^T$)

We give here some particular signals in $F$ as well as some properties used in section 3 (traffic model).

- The gain signal $\gamma^p$:

$$\gamma^p(t) = \begin{cases} p \text{ if } t = 0 \\ +\infty \text{ for } t > 0 \end{cases}$$

- The shift signal $\delta^T$:

$$\delta^T(t) = \begin{cases} 0 \text{ if } t \leq T \\ +\infty \text{ otherwise} \end{cases}$$

It is then easy to obtain the following signals:

- The signal $\gamma^p * \delta^T$:
$$(\gamma^p * \delta^T)(t) = \begin{cases} p \text{ if } t \in [0, T] \\ +\infty \text{ if } t > T \end{cases}$$

- The signal $(\gamma^p * \delta^T)^*$:
$$(\gamma^p * \delta^T)^*(t) = \begin{cases} 0 \text{ if } t \leq 0 \\ kp \text{ for } kT \leq t < (k+1) > T, \quad k \in N \end{cases}$$

Moreover, we can easily check (see Figure 7) that
$$(\gamma^p * \delta^T)^*(t) \geq \frac{p}{T} t, \quad \forall t.$$

We explain here the convolution of the signals $\gamma^p$ and $\delta^T$ with a signal $f$.

- $(\gamma^p f)(t) := (\gamma^p * f)(t) = f(t) + p, \forall t.$

- $(\delta^T f)(t) := (\delta^T * f)(t) = f(t - T), \forall t.$

- $(\gamma^p \delta^T f)(t) := (\gamma^p * \delta^T * f)(t) = f(t - T) + p, \forall t.$

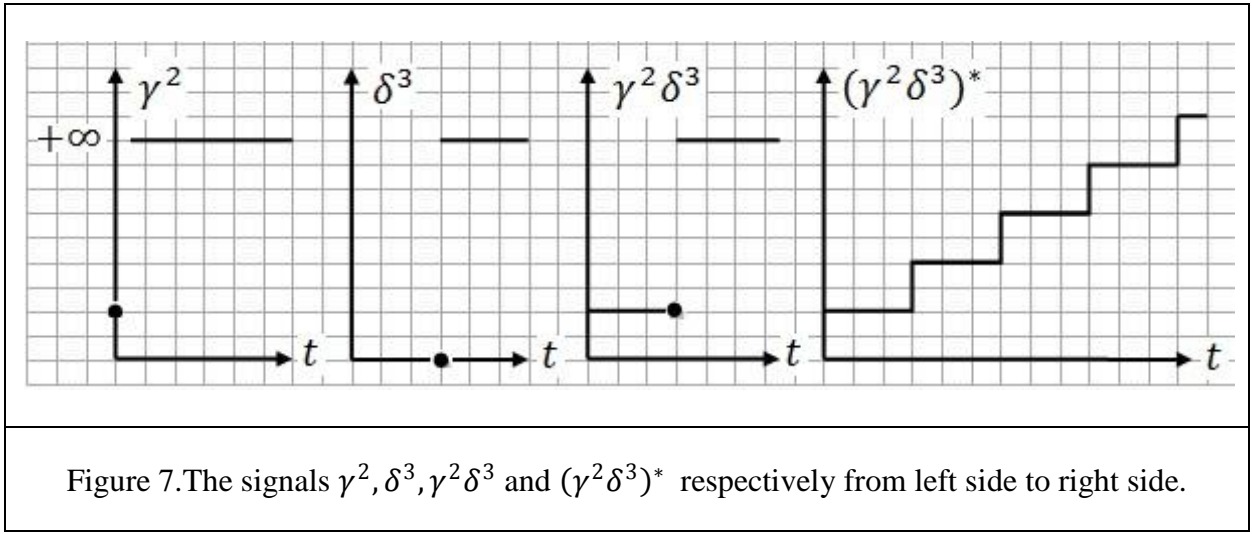

Figure 7. The signals $\gamma^2, \delta^3, \gamma^2\delta^3$ and $(\gamma^2\delta^3)^*$ respectively from left side to right side.

We recall the following additional properties; see (Baccelli, Cohen, Olsder, & Quadrat, 1992) and/or (Le Boudec & Thiran, 2001).

- (P1): $\forall f \in F, f^* \leq f * f^* \leq f.$
- (P2): $\forall f \in F, f(0) = 0 \Rightarrow f * f^* = f^*.$
- (P3): $\forall f, g \in F, (f \oplus g)^* = f^* * g^*.$
- (P4): $\forall f, g \in F, f \geq 0, g \geq 0 \Rightarrow [(f * g) * \gamma^{-a}]^+ = f * [g * \gamma^{-a}]^+.$
- (P5): $\forall f \in F, e \oplus f * f^* = f^*.$

## Appendix B (minimum service couple)

We clarify here the three bounds of network calculus in the case where we have a minimum service couple $(\beta, \lambda)$ instead of a minimum service curve $\beta$. That is to say that we have $Y \geq \beta * U \oplus \lambda$ instead of $Y \geq \beta * U$.

- The maximum backlog:
$$B(t) = U(t) - Y(t) \leq U(t) - \min\left\{\min_{0 \leq s \leq t}(U(t-s) + \beta(s)), \lambda(t)\right\}$$
$$\leq \max\left\{\max_{0 \leq s \leq t}(\alpha(s) - \beta(s)), \alpha(t) - \lambda(t)\right\}$$
$$\leq \max\left\{\max_{z \geq 0} \max_{0 \leq s \leq z}(\alpha(s) - \beta(s)), \max_{z \geq 0}(\alpha(z) - \lambda(z))\right\}$$
$$\leq \max\left\{\max_{s \geq 0}(\alpha(s) - \beta(s)), \max_{s \geq 0}(\alpha(s) - \lambda(s))\right\}$$
$$\leq \max_{s \geq 0}\{\alpha(s) - (\beta \oplus \lambda)(s)\}.$$

- The maximum delay:
  Let $t \geq 0$. Let $\tau \geq 0$ such that $\tau < d(t) := \inf\{h \geq 0, U(t) \leq Y(t+h)\}$. Then we have $U(t) > Y(t+\tau)$. We have
  $$Y(t+\tau) \geq \min\left\{\min_{0 \leq s \leq t+\tau}\{U(t+\tau-s) + \beta(s)\}, \lambda(t+\tau)\right\}.$$
  That is
  $$\exists s_0, 0 \leq s_0 \leq t+\tau, \quad Y(t+\tau) \geq \min\{U(t+\tau-s_0) + \beta(s_0), \lambda(t+\tau)\}.$$
  Then
  $$\exists s_0, 0 \leq s_0 \leq t+\tau, \quad U(t) \geq \min\{U(t+\tau-s_0) + \beta(s_0), \lambda(t+\tau)\}.$$
  Therefore
    - If $U(t) \geq U(t+\tau-s_0) + \beta(s_0)$ then $s_0 > \tau$ and thus
      $$\alpha(s_0 - \tau) \geq U(t) - U(t+\tau-s_0) > \beta(s_0)$$
      Hence
      $$\tau < \sup_{t \geq 0} \inf\{h \geq 0, \alpha(t) \leq \beta(t+h)\} \leq \sup_{t \geq 0} \inf\{h \geq 0, \alpha(t) \leq (\beta \oplus \lambda)(t+h)\}.$$
    - If $U(t) \geq \lambda(t+\tau)$ then since $U(0) = 0$, we have $\alpha(t) \geq U(t) \geq \lambda(t+\tau)$.
      Hence
      $$\tau < \sup_{t \geq 0} \inf\{h \geq 0, \alpha(t) \leq \lambda(t+h)\} \leq \sup_{t \geq 0} \inf\{h \geq 0, \alpha(t) \leq (\beta \oplus \lambda)(t+h)\}.$$

- The departure flow:
$$Y(t) - Y(s) \leq U(t) - \max\left\{\min_{0 \leq z \leq s}(U(s-z) + \beta(z)), \lambda(s)\right\}$$
$$\leq \max\left\{\max_{0 \leq z \leq s}\{U(t) - U(s-z) - \beta(z)\}, U(t) - \lambda(s)\right\}$$
$$\leq \max\left\{\max_{0 \leq z \leq s}\{\alpha((t-s) + z) - \beta(z)\}, \alpha(t) - \lambda(s)\right\}$$
$$\leq \max\left\{\max_{0 \leq z \leq s}\{\alpha((t-s) + z) - \beta(z)\}, \max_{0 \leq z \leq s}\{\alpha((t-s) + z) - \lambda(z)\}\right\}$$
$$\leq \max_{0 \leq z \leq s}\{\alpha((t-s) + z) - \min(\beta, \lambda)(z)\} \leq (\alpha \oslash (\beta \oplus \lambda))(t-s).$$